\documentclass[12pt]{article}
\usepackage{amsmath, amssymb, amsthm}
\usepackage{geometry}
\usepackage{hyperref}
\usepackage{tcolorbox}

\title{The Inverse Function Fallacy:\\ On Sign Determination and Forgotten Fundamentals}
\author{Meliksah Yorulmazlar \\
Department of Computer Science,\\
Rensselaer Polytechnic Institute \\
110 8th Street, Troy, NY 12180, USA \\
\texttt{yorulk@rpi.edu}}
\date{17 September 2025}

\begin{document}
\maketitle

\begin{abstract}
Although inverse functions are introduced early in algebra, many students remain unaware that an inverse expression may legitimately involve a negative root. Instead, they default to assuming a positive root, overlooking the role of domain restrictions in determining the correct solution. This paper identifies this misconception as the \textit{inverse function fallacy} and introduces a systematic approach---the \textit{Inverse Function Point} method---that establishes sign determination to a single domain-based reference point. In a study of 69 STEM students at Rensselaer Polytechnic Institute, only 19\% solved a sign-determination problem correctly on their first attempt. Those students were not re-tested. Of the remaining 56 who answered incorrectly, I was able to re-test 40 after teaching the proposed method. In this subgroup, accuracy rose to 80\%. These results highlight both the fragility of assumed mathematical knowledge and the potential of simple, intuitive procedures to reinforce conceptual understanding.
\end{abstract}

\begin{quote}
``It ain't what you don't know that gets you into trouble. It's what you know for sure that just ain't so.'' --- Mark Twain
\end{quote}

\section{Introduction}
While revisiting my A-level mathematics coursebook \cite{pemberton2018}, I encountered inverse function problems that I struggled to solve intuitively. In particular, I could not determine which sign of the root to choose when the algebraic manipulation yielded both a positive and negative solution. Consulting the solution manual provided no systematic method, only case-by-case reasoning.

\medskip
This frustration led me to reflect more deeply on the structure of inverse functions. The geometric interpretation is that an inverse function corresponds to reflecting a function across the line $y = x$. This perspective, combined with domain restrictions, makes it possible to determine the correct root logically rather than by guesswork or trial-and-error.

\medskip
Motivated by this insight, I conducted a study at Rensselaer Polytechnic Institute (RPI). At first, I wanted to confirm whether STEM students actually shared this misunderstanding. Once the fallacy was observed, I then tested whether the proposed method could help resolve it.

\section{Background}
Inverse functions are typically introduced in secondary school, with domain restrictions emphasized to ensure one-to-one correspondence. Yet even at the university level, many students fail to recognize that an inverse function may legitimately involve a negative root. Instead, they default to assuming a positive root, overlooking the role of domain restrictions in determining the correct root.  

\medskip
This study investigates the prevalence of this misconception and introduces a method designed both to make students aware of the possibility of negative roots and to provide a straightforward procedure for identifying the correct one. 

\medskip
Prior research in mathematics education has highlighted persistent misconceptions around functions and inverses. Vinner and Dreyfus \cite{vinner1989} described the gap between students’ \textit{concept images} (personal mental representations) and formal \textit{concept definitions}, showing how this mismatch often leads to systematic errors. This study adds to that literature with empirical data from university-level STEM students.

\section{Methodology}
I designed worksheets containing five problems: three requiring sign determination (two with negative roots, one with a positive root), one involving solving for $x$, and one requiring use of completing the square. Each worksheet also included a single worked example to remind participants of the basic process of finding an inverse function, though this example did not feature the proposed method. Those who agreed were given a worksheet and completed it on site. The problems were distributed informally during RPI’s summer tutoring sessions, which are held in the student union. During these sessions, students typically work independently on assignments or seek help from tutors. When a student appeared available, I explained that the worksheet was part of a short research study on inverse functions and invited them to participate. Those who agreed were given a worksheet.

\medskip
Participants ($n=70$) were self-selected volunteers, predominantly STEM majors, ranging from junior undergraduates to PhD students. One non-STEM participant was excluded for homogeneity, leaving $n=69$. For each participant, I recorded both their program level (UG, MS, MENG or PhD) and whether they were a mathematics major.

\medskip
No time limits or restrictions were imposed; participants were free to use calculators and to discuss with peers. In principle, the type of problems posed—algebraic manipulations and inverse function exercises—could have been checked quickly using AI-based tools, since students had access to phones or laptops. However, the informal tutoring setting and the low-stakes nature of the activity gave them little incentive to do so. As a result, the low scores observed among participants can be taken as a genuine indicator of conceptual misunderstanding, rather than artifacts of testing conditions, time pressure, or reliance on external assistance.

\section{Results}
 Of the 69 participants included in the study, only 13 (19\%) correctly determined the appropriate sign when solving inverse function problems on their first attempt. These students were not re-tested, having already demonstrated understanding. The vast majority either assumed the positive root by default or were unaware that a negative root was even possible. This indicates that the misconception is not simply about applying the wrong rule, but about failing to recognize the existence of sign ambiguity in the first place.  

\medskip
Among the correct solvers, the most common strategy was to sketch the original function, reflect it across $y = x$, and then restrict the domain to verify the correct root. One undergraduate even reported relying on Desmos to graph the function in order to intuit which sign to choose. While these approaches were valid, they proved highly time-consuming and impractical for routine problem solving, especially on timed assessments. This reliance on graphing underscored the need for a method that can determine the correct sign algebraically, without requiring a sketch or computational tool. Notably, even some mathematics PhD students failed to arrive at the correct solution,
highlighting the persistence of this conceptual gap across levels of expertise.

\medskip
Additional unexpected observations also emerged. Several undergraduates reported unfamiliarity with the notation $gf(x)$, having only encountered $g \cdot f(x)$ or $g(f(x))$ in prior coursework. Others struggled with completing the square, despite its status as a standard high school algebra technique. These difficulties suggest that the inverse function fallacy is connected to broader fragility in mathematical foundations.  

\medskip
Overall, the initial results confirmed both the prevalence and the persistence of the inverse function fallacy. The data underscored the need for a streamlined method, one that could help students bypass the reliance on graphing and avoid systematic sign errors.

\medskip
 These findings highlight the gap between procedural skill and conceptual understanding, providing the rationale for the method introduced in the next section.
\section{Proposed Method}
A common misconception about inverse functions is that students assume the inverse always carries a positive root. In reality, inverse functions may require either a positive $(+)$ or negative $(-)$ root, depending on the domain of the original function. The method presented here introduces the concept of the \textit{Inverse Function Point}, which makes the correct sign evident without graphing.

\medskip

\textbf{Example:} Consider $f(x) = x^2 - 2, \; x \leq 0$. Find $f^{-1}(x)$.

\begin{enumerate}
    \item \textbf{Select a domain point.}  
    Since the domain is $x \leq 0$, we may choose any valid input, such as $x = -2$.  
    Then $f(-2) = (-2)^2 - 2 = 2$, which gives the point $(-2,2)$ on $f(x)$.

    \item \textbf{Form the Inverse Function Point.}  
    By reflecting across $y = x$, the coordinates swap, giving the \textit{Inverse Function Point} $(2,-2)$.  
    This point must lie on the inverse, and will later confirm which root is valid.

    \item \textbf{Solve algebraically for the inverse.}  
    Starting from $y = x^2 - 2$:  
    \[
    x = y^2 - 2 \quad \Rightarrow \quad y = \pm \sqrt{x+2}.
    \]

    \item \textbf{Confirm using the Inverse Function Point.}  
    Substitute the $x$-value of the Inverse Function Point ($x=2$) into the candidate inverses:  
    \[
    y = \pm \sqrt{2+2} = \pm 2.
    \]  
    Only the negative root matches the $y$-value of the Inverse Function Point $(2,-2)$.  
    Therefore, the correct inverse is:
    \[
    f^{-1}(x) = -\sqrt{x+2}, \quad x \geq -2.
    \]
\end{enumerate}

\medskip

By starting the process with a single Inverse Function Point, students can determine the correct sign logically. This removes the guesswork, reinforces the reflection principle, and avoids the need for graphing for complicated problems.

\medskip
\begin{tcolorbox}[colback=blue!5,colframe=blue!40!black,title=Rule: Inverse Function Point Method]
To determine the correct sign of an inverse function:
\begin{enumerate}
    \item Select a point from the restricted domain of the original function.  
    \item Reflect it across $y=x$ to form the \textit{Inverse Function Point}.  
    \item Solve algebraically for the inverse (which will usually produce $\pm$ roots).  
    \item Substitute the $x$-value of the Inverse Function Point; the correct root is the one whose $y$-value matches.  
\end{enumerate}
\end{tcolorbox}

\section{Follow-up Study}
To test the effectiveness of the proposed method, I prepared a second worksheet that contained a one-page explanation of the sign determination procedure followed by four practice problems. The explanation emphasized the role of domain restrictions and the reflection principle across $y=x$, guiding students to identify whether the positive or negative root was appropriate. 

\medskip
The 13 students who solved the initial problem correctly had already demonstrated mastery and were therefore excluded from the follow-up. Of the 56 remaining, 40 were re-tested under similar conditions. Participants first read a one-page explanation of the method and then attempted four new problems.

\medskip
Results showed a substantial improvement: 32 out of 40 participants (80\%) successfully solved the problems on their second attempt. Several students noted in optional comments that the reflection-based approach clarified why a negative root can be valid, and that testing an inverse point to confirm the sign was a new idea.

\medskip
These findings indicate that the method not only improves accuracy but also promotes a deeper conceptual understanding. By revealing to students that a negative root is possible and providing a systematic way to decide between signs, the method addresses both layers of misunderstanding observed in the initial study.

\section{Discussion}
These results confirm the persistence of the inverse function fallacy even among advanced learners. The method presented offers a systematic shortcut, analogous to algebraic rules replacing graphing in differentiation. Broader implications include:
\begin{itemize}
    \item Many students were not only uncertain about which sign to select, but were unaware that a negative root could exist at all, revealing a deeper conceptual gap.
    \item The fragility of foundational knowledge (e.g., completing the square). 
    \item The pedagogical importance of reinforcing intuitive principles like reflection across $y=x$. 
    \item The potential to impact tens of millions of students worldwide, given that inverse functions are a standard component of secondary mathematics curricula.

\end{itemize}

\medskip
The inverse function fallacy illustrates Twain’s observation: students were not hindered by lack of exposure to inverses, but by the false certainty that the inverse always carries a positive root.

\section{Conclusion}
This work highlighted what I call the \textit{inverse function fallacy} and presented a practical method for resolving sign ambiguities using domain-based reflection. A key finding was that many students were not only uncertain about which sign to choose, but were unaware that a negative root could exist at all. Experimental evidence with RPI students showed significant improvement after being introduced to the method. With this new insight, students can better learn and intuit the concept of inverse functions, while teachers gain a more systematic and intuitive way to present it. Revisiting and reinforcing such fundamentals is essential for strengthening mathematical intuition and supporting advanced problem-solving.

\section*{Acknowledgments}
I am grateful to all the students at Rensselaer Polytechnic Institute who participated in this study and contributed their time and insights. Their willingness to engage with foundational mathematics problems provided the basis for this research.

\medskip

I would also like to extend special thanks to \textbf{Mehtap Agirsoy}, PhD student in Mechanical Engineering, and \textbf{Eliane Traldi}, Lecturer in Mathematics, for their thoughtful feedback, encouragement, and support throughout this project.

\medskip
I would also like to thank the many books I have read, which have taught me to look beyond what is given and to question the world around me. They have instilled in me the habit of seeking deeper understanding, a habit that guided me in uncovering the ideas explored in this work.

\section*{Author Statement}
A large language model was used for limited language editing and clarification of wording.
All data collection, analysis, interpretation, and conclusions were performed by the author.
The author takes full responsibility for the content of the manuscript.

\section*{Ethical Considerations}
The study involved voluntary participation by adult students in an informal, non-instructional setting. No identifying information was collected, participation was optional, and no incentives or academic credit were offered. The activity posed no more than minimal risk to participants. Under institutional guidelines, this type of informal observational data collection did not require formal ethics board approval.

\end{document}